\documentclass[12pt]{article}
\usepackage{epsfig}
\usepackage{amsmath,amsfonts}
\usepackage{natbib,apalike}

\def\e{\hbox{E}}
\def\cov{\hbox{Cov}}
\def\corr{\hbox{Corr}}
\def\var{\hbox{Var}}

\def\max{\hbox{max}}
\def\key{{\cal K}}
\def\Beta{\hbox{Beta}}

\def\bx{\textbf{x}}
\def\by{\textbf{y}}
\def\bu{\textbf{u}}

\def\bX{\textbf{X}}
\def\bK{\textbf{K}}
\newcommand{\bmu}{\mbox{\boldmath $\mu $}}
\def\bSigma{{\bf \Sigma}}
\def\bvn{\hbox{BVN}}
\def\mvn{\hbox{MVN}}
\def\b0{{\bf 0}}

\begin{document}

\title{Inference for Quantile Measures of Kurtosis,
Peakedness and Tail-weight}

\author{Robert G. Staudte\\
La Trobe University, Melbourne, Australia\thanks{Emeritus Professor
Robert G. Staudte, Department of Mathematics and Statistics, La~Trobe University, Melbourne, Vic. 3086 Australia, {\em r.staudte@latrobe.edu.au}}}
\date{23 July, 2014}
\maketitle

\clearpage
\newpage

\abstract{Many measures of peakedness, heavy-tailedness and kurtosis have been proposed in the literature, mainly because kurtosis, as originally defined, is a complex combination of the other two concepts. Insight into all three concepts can be gained by studying Ruppert's ratios of interquantile  ranges.
They are not only monotone in Horn's measure of peakedness when applied to the central portion of the
population, but also monotone in the practical tail-index
of Morgenthaler and Tukey, when applied to the tails.
Distribution-free confidence intervals are found for Ruppert's ratios, and sample sizes required to obtain such intervals for a pre-specified relative width and level are provided. In addition, the empirical power of distribution-free tests for peakedness and bimodality are found for symmetric beta families and mixtures
of $t$ distributions.  An R script that computes the
confidence intervals is provided in online supplementary material.
}

\vspace{3cm}


\clearpage
\newpage

\section{INTRODUCTION}\label{intro}

\subsection{Background and summary}

The meaning of kurtosis has long puzzled statisticians, ever since the standardized fourth moment definition was introduced by \cite{kpear-1905} to help describe departures from normality.   A century elapsed before its asymptotic distribution was derived by \cite{pewsey-2005}, although
its sister sample skewness result was obtained much earlier in \cite{gupta-1967}.
In the meantime, numerous other measures of kurtosis have been proposed and dissected, but again with almost no accompanying inferential methods.

There are three themes pervading research into kurtosis measures.  Firstly, kurtosis as originally
conceived is a location, scale and sign-invariant measure of shape that somehow measures both peakedness and tail-weight. Contributions by many authors to this theme are thoroughly described by \cite{bala-1988}. Secondly, an increase in a kurtosis measure should quantify movement of mass from the tails to the center of the distribution,  with substantive contributions from \cite{zwet-1964}, \cite{oja-1981} and  \cite{bala-1988, bala-1990}.
The third theme is that quantile-based measures are preferable to moment-based measures: they are always defined and are robust in that they have bounded influence functions and positive breakdown points.
Contributions of this type include \cite{groen-1984}, \cite{rupp-1987}, \cite{moors-1988}, \cite{groen-1998}
and \cite{kotz-2009}. Also of interest are  the maximum-bias curves for interquantile ranges studied by \cite{croux-2001}, the robust kurtosis measures of
  \cite{seier-2003},  and the  $L$-moment kurtosis measures of \cite{with-2011}.

 Recently \cite{JRP-2011} studied ratios of linear combinations of interquantile ranges, and showed that
they possessed the surprising property of invariance to skewness-inducing transformations. The simplest
measures of this type, ratios of two interquantile ranges, were introduced by
\cite{rupp-1987}, who compared their influence functions and order-preserving properties with other measures of kurtosis. Despite their simplicity, they provide a basis for studying the peakedness and tail-weight
properties of distributions, separately or jointly.

As explained further in Section~\ref{peaktail}, these simple ratios  measure peakedness when applied to
the center of a distribution, and they measure tail-weight when applied  to the remaining (tails) portion.
This idea is already exploited by \cite{schmid-2003}, who found tests for normality based on these ratios
of ranges. In Section~\ref{peakedness} we extend the peakedness measure of \cite{horn-1983} so that it can
detect bimodality, and show that the \cite{rupp-1987} kurtosis, when applied to the central portion of the
distribution, continues to be approximately monotone in it. We further show in Section~\ref{tailweight}
that, when applied to the tails portion, the Ruppert kurtosis is monotone in the index of tail-weight of \cite{M-T-2000}.

In Section~\ref{Funknown} we briefly describe inference for the ratio of interquantile ranges when the
underlying location-family is known; it is based on a variance stabilizing transformation (VST)  which requires three constants, each depending on the family through the sparsity index of \cite{tukey-1965}. By estimating these constants (nuisance parameters), which requires density estimates at four quantiles,
one obtains distribution-free confidence intervals for the ratio of interquantile
ranges.  These intervals are evaluated by simulation studies for  coverage and widths in Section~\ref{simulations}. The coverage for  90\% or 95\% confidence intervals is accurate provided that the sample size is at least 400. The empirical power of the Ruppert measures for detecting peakedness and/or bimodality is also found for the symmetric Beta models and mixtures of $t$ distributions.
 A summary and further research problems are outlined in Section~\ref{summary}.

\subsection{Preliminary definitions and concepts}

 For any strictly increasing distribution function $F$ and $0<t<1$ let $x_t=G(t)\equiv F^{-1}(t)$ denote the $t$th quantile. For $0<t<0.5$ denote the $t$th {\em interquantile range} of $F$ by $R_t=R_t(F)=x_{1-t}-x_t $. Then for $0<p<r< 1/2$ \cite{rupp-1987} defined a measure of kurtosis by $\kappa _{p,r}=R_p/R_r.$ (Our notation differs from his: our $\kappa _{p,r}$ is his $R_{r,p}.$)
 These measures are clearly sign, location and scale invariant. Our choice of $(p,r)$ is guided by a desire to have a quantile measure which agrees at the normal model
    with the classical
   moment-based definition of kurtosis $\alpha _4(F)=\mu _4/\mu _2^2$, where $\mu _k=\e _F[(X-\e [X])^k]$, is the
   $k$th moment about the mean $\e _F[X],$ $k=2,3,\dots .$  The normal model $F=\Phi$
    has $\alpha _4(\Phi )=3$.
   In the case of symmetric $F$, $\kappa _{p,r}=x_p/x_r,$ so to have $\kappa _{p,r}=3$ for the normal distribution, we need to have $p=p(r)= \Phi (3\Phi ^{-1}(r))$. Some examples are given in Table~\ref{table1}.
     Further, we want to be able to carry out tests and
   find confidence intervals for $\kappa _{p,r},$ and to provide some protection against outliers by
   choice of $(p,r)$.
\begin{table}[h!]
\caption{\label{table1}\em \footnotesize Examples of the kurtosis coefficient $\kappa _{p,r}=R_{p}/R_{r}$, for various models and four choices of\; $(p(r),r)$, with $r=0.3, 0.333, 0.35,0.4$ and $p(r)=\Phi(3\Phi ^{-1}(r))$. Also shown is the classical kurtosis $\alpha _4(F)$.}
\begin{small}
\begin{center}
\begin{tabular}{lrrrrr}
 \hline
 $\qquad F$ & $\quad \alpha _4(F)$ &  $r=0.3$ & $r=1/3$ & $r=0.35$ &  $r=0.4$ \\
  \hline
1. Beta$(1/2,1/2)$  & 1.50 & 1.673  & 1.906 &  2.038   &  2.470\\
2. Uniform          & 1.80 & 2.211  & 2.411 &  2.508   &  2.764\\
3. Beta$(2,2)$      & 2.14 & 2.588  & 2.709 &  2.762   &  2.892\\
4. Normal           & 3.00 & 3.000  & 3.000 &  3.000   &  3.000\\
5. Logistic         & 4.20 & 3.294  & 3.200 &  3.160   &  3.070\\
6. Student-$t_5$    &  9.00& 3.399  & 3.260 &  3.205   &  3.086\\
7. Student-$t_4$    & $-$  & 3.523  & 3.337 &  3.265   &  3.110\\
8. Student-$t_2$    & $-$  & 4.340  & 3.820 &  3.631   &  3.250\\
9. Laplace          & 6.00 & 4.223  & 4.016 &  3.913   &  3.606\\
10. Cauchy          & $-$  & 7.492  & 5.438 &  4.787   &  3.635\\
\noalign{\smallskip}
11. Beta$(2,1)$     & 2.40 &  2.527 &   2.661 &   2.722 &  2.872   \\
12. $\chi ^2_5$     & 5.40 &  3.088 &   3.060 &   3.048 &  3.021   \\
13. $\chi ^2_3$     & 7.00 &  3.167 &   3.113 &   3.091 &  3.039   \\
14. $\chi ^2_2$     & 9.00 &  3.293 &   3.200 &   3.161 &  3.070   \\
15. $\chi ^2_1$     &15.00 &  3.881 &   3.625 &   3.511 &  3.232   \\
16. Log-normal      &113.94&  4.205 &   3.789 &   3.624 &  3.262   \\
17. Skew-$t_{2,2}$  &$-$   &  4.340 &   3.820 &   3.631 &  3.250   \\
18. Pareto(2)       &$-$   &  4.961 &   4.216 &   3.941 &  3.377   \\
19. Skew-$t_{2,1}$  & $-$  &  7.492 &   5.438 &   4.787 &  3.635   \\
20. Skew-$t_{2,1/2}$& $-$  & 30.452 &  14.033 &  10.189 &  4.984   \\
   \end{tabular}
  \end{center}
  \end{small}
     \end{table}

The models in Table~\ref{table1}  are labeled with standard notation in \cite{J-K-B-1994,J-K-B-1995}, but two cases require clarification: the Pareto distribution with shape parameter $a=2$ has  distribution function given by $F(x)=1-1/x^2$ for $x\geq 1.$ The class of \lq skewed-$t$\rq\  distributions introduced by \cite{RJP-2011} are denoted $t_{\epsilon ,\nu}$ where $\epsilon $ is a real skewness parameter and $\nu >0$ is the degrees of freedom. If $X\sim t_\nu $, then
$Y=\sinh(\sinh ^{-1}(X)+\epsilon )\sim t_{\epsilon ,\nu}.$ Clearly $t_{0 ,\nu}=t_\nu $ and $t_{\epsilon ,+\infty }$ is the skewed normal model, while $t_{\epsilon,1}$ is the skewed Cauchy model.  A nice property of these distributions
is that ratios of linear combinations of  interquantile ranges  are not dependent on the skewness
parameter $\epsilon $, see \cite{JRP-2011}.  However, as shown in Section~\ref{simulations}, $\kappa _{p,r}(t_{\epsilon ,1})$ is much more difficult
to estimate for $\epsilon =2$ than $\epsilon =0$.

The second column of Table~\ref{table1} gives values of the classical moment kurtosis for the Models in
Column 1.  The remaining columns give values of $\kappa _{p(r),r}$ for $r=0.3, 1/3 , 0.35,0.4$.
Note that the  kurtosis $\kappa _{p(r),r}$ becomes more discriminating as $r$ gets smaller; however even for $r=0.25$, the value of $p(r)$ is 0.0215, so $r<0.3$ is excluded to guarantee resistance to 5\% of outliers.
  As $r\to 0.5$, $\kappa _{p(r),r}\to 3$, by L'Hospital's rule. Therefore larger values of $r\geq 0.4$ are less informative. Within the range  $0.3 \leq r \leq 0.4$ we  decided to focus on $r=1/3$ because then the ordering
  of  $\kappa _{1/3}=\kappa_{p(1/3),1/3}(F)$ for the various models $F$ in Table~\ref{table1}
   is roughly consistent with that of $\alpha _4(F)$, as well as agreeing exactly at $F=\Phi $. Further, it is easy to remember that because $p(1/3)\approx 0.1$, one is comparing the range of the middle 4/5 of the
   population with the range of the middle 1/3.

\section{PEAKEDNESS AND TAILWEIGHT}\label{peaktail}

Throughout this section fix $0 < p < q <r < 0.5.$  The \lq central\rq\ portion of the distribution
of $F$ is that lying between $x_q$ and $x_{1-q}$ while the \lq tail\rq\ portion is that lying outside these
quantiles. We will show that applying the  kurtosis measure of \cite{rupp-1987} to the center of the distribution leads to a peakedness measure, while applying it to the tails portion leads to a tail-weight measure.
To this end, define the (central) \textit{quantile peakedness} by $\pi_{q,r}=R_q/R_r$, for $q<r<0.5.$ Define the \textit{quantile tail-weight} by $\tau_{p,q}=R_p/R_q$, for $0<p<q$.  Trivially, the product is the \lq kurtosis\rq\  measure $\kappa _{p,r}=\tau_{p,q}\,\pi_{q,r}=R_p/R_r$ for the distribution $F$.
All three  measures satisfy the kurtosis convexity criterion of \cite{zwet-1964} and that of
\cite{lawr-1975}, see \cite[Theorem 2]{rupp-1987}. And each has the simplest form of a skewness invariant
kurtosis measures \cite[Sec. 2.1]{JRP-2011}.

\cite{schmid-2003} carried out tests for peakedness, tail-weight and leptokurtosis based on sample versions
of $\pi_{q,r}$, $\tau_{p,q}$ and  $\kappa _{p,r}$, respectively, for the case of $p=1/40$, $q=1/8$ and
$r=1/4$ (our notation).  We prefer larger
values because their choice of $r=1/4$ means that the central half of the data are ignored in assessing peakedness. Further, their choice of $p=1/40$ means that the breakdown point of the procedure is only 1/40.
In any case, their emphasis is on testing while ours is on confidence intervals so the results to follow
can be seen as complementary to theirs.

\subsection{Peakedness Measures}\label{peakedness}
A justification for calling  $\pi_{q,r}=R_q/R_r$ \textit{quantile peakedness} for symmetric unimodal distributions is already given by
 \cite{rupp-1987}.  He showed that for $r$ less than, but near 0.5, $\pi_{q,r}$ is approximately monotone increasing  in the peakedness measure  of \cite{horn-1983}. However, Horn only considered symmetric densities $f$  that were unimodal. Next we extend his measure of central peakedness to one that distinguishes bimodality and show that $\pi_{q,r}$ is still approximately monotone increasing in this extended version.

\begin{figure}[t!]
\centering
\includegraphics[scale=.8]{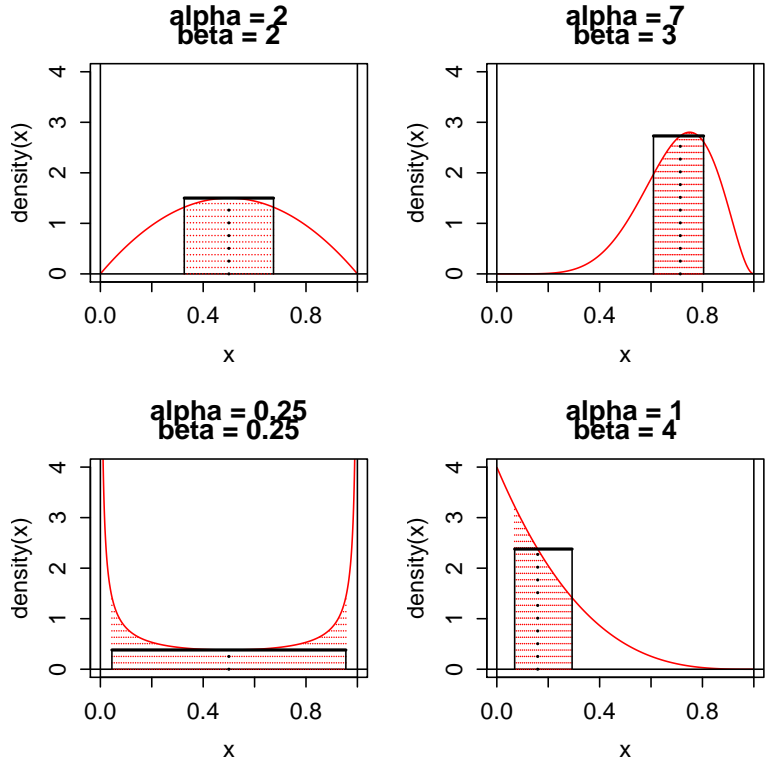}
\caption{\footnotesize In these plots $q=1/4$. Four $\Beta (\alpha ,\beta $) densities are shown, with parameters listed above each plot. The areas lying under the densities and over the interval $[x_q,x_{1-q}]$ are shaded and have areas equal to $1-2q=1/2.$ These are to be compared to the areas
of the rectangles $A_q=f(x_{0.5}; \alpha ,\beta)\; R_{1/4}(\alpha ,\beta )$. The respective medians are marked by the vertical dotted  lines. See text for more details.
\label{fig1}}
\end{figure}

\begin{figure}[t!]
\centering
\includegraphics[scale=.8]{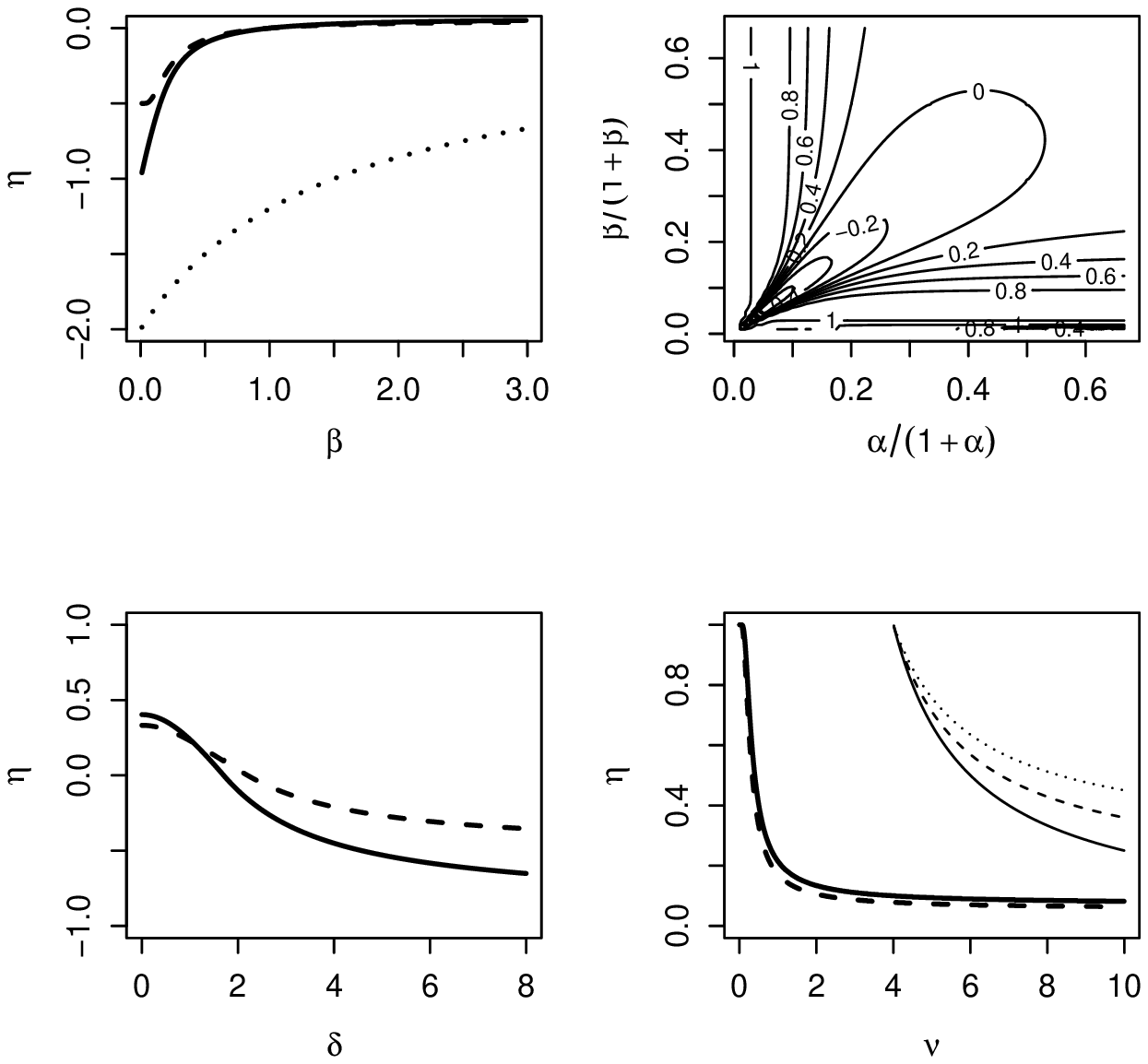}
\caption{\footnotesize In these plots $q=1/4$ and $r=3/8.$ In the top left plot is shown the graph of $\eta _q$ defined in (\ref{hornext}) as a function of $\beta $ for the $\Beta (\beta ,\beta )$ model as a thick solid line.
Its approximation $\hat \eta _{q,r}$  is plotted as a thick dashed line. The dotted line
shows the graph of $\alpha _4-3$.  In the upper right plot are shown contours of $\eta _q$ for the
$\Beta (\alpha ,\beta )$ model. Note that it takes on negative values within the contour marked 0.
The bottom left plot again shows $\eta _q$ and its approximation $\hat \eta _{q,r}$ for a 50:50 mixture of two $t_{1/2}$ distributions that are distance $\delta $ apart. The bimodality is detected in that $\eta _q<0$ for $\delta >1.5.$  The approximation of $\eta _q$ by $\hat \eta _{q,r}$ improves as $r$ moves closer to 0.5. The bottom right plot shows $\eta _q$ and its approximation as functions of $\nu $ for
the skew-$t_{\epsilon ,\nu}$ distributions; it does not depend on $\epsilon .$  Also shown are values of
$(\alpha _4-3)/\alpha _4$ for $\epsilon =0, $ 0.25, and 1, respectively in thin solid, dashed and dotted lines.
\label{fig2}}
\end{figure}

\subsubsection*{A simple extension of Horn's measure of peakedness}\label{horn}

 \cite{horn-1983} considered densities such as that depicted in the upper left plot of Figure~\ref{fig1}.
 Consider the rectangle with base $[x_{q}, x_{1-q}]$  and height $f(x_{0.5})$ which has area $A_q=f(x_{0.5})R_q$; then Horn's measure (our notation) is based on the  ratio $(1-2q)/A_q$,
 which for symmetric unimodal $f$ is the proportion of the area
of the rectangle which lies under the density. Clearly
this ratio, which lies between 0 and 1, will be {\em smaller} with more peakedness. To make it increasing
in peakedness,  \cite{horn-1983} defined $\eta _q=1-(1-2q)/A_q$, which still varies from 0 to 1, but now with {\em larger} values indicating more peakedness.

 At the other extreme, the bottom left plot in Figure~\ref{fig1} indicates that symmetric U-shaped distributions with minimum at the median will have  the ratio $(1-2q)/A_q >1.$ These observations motivate
 a measure of (central) peakedness defined by
\begin{eqnarray}\label{hornext}
\eta  _q
 &= & \left\{
              \begin{array}{ll}
   -1+A_q/(1-2q)       & \hbox {\quad for $A_q \leq (1-2q) $\, ;} \\
        +1-(1-2q)/A_q,     & \hbox {\quad for $ (1-2q)\leq A_q $\,.}
              \end{array}
            \right.
 \end{eqnarray}
 This $\eta  _q$  agrees with Horn's definition for symmetric unimodal $f$, but can be applied to arbitrary $f$, even if $f(x_{0.5})=0$ or $+\infty $. It lies in $[-1,1]$, takes on negative values for symmetric U-shaped models, and equals 0 for the uniform distribution.

Some examples of $\eta  _q$ for $q=1/4$ are shown in Figure~\ref{fig2}, where its graph is plotted as a solid line for the $\Beta(\beta,\beta )$, $\beta >0$ family; a 50:50 mixture of two Student-$t_{1/2}$ models, one of which is shifted by $\delta >0$, and the skew-$t_{\epsilon ,\nu }$ families
for $0<\nu <10$ and selected values of $\epsilon $. Also shown is a contour plot of $\eta  _q$ for the
 $\Beta(\alpha,\beta )$  family, $\alpha >0,\beta >0$.    These plots confirm that $\eta _q$ can detect bimodality as well as peakedness.

\subsubsection*{The measure $\pi _{q,r}$
is monotone increasing in $\eta  _q$\,.}
 An approximation to $\eta  _q$ can be obtained as  in \cite{rupp-1987}: for small $\epsilon >0$
 one has the finite difference approximation $f(x_{0.5})\approx 2\epsilon/\{ x_{0.5+\epsilon}-x_{0.5-\epsilon } \}.$  Thus for $0<q<r< 0.5$ and $r$ near 0.5, say $\epsilon =0.5-r$,
\begin{equation}\label{ratioapprox}
    \frac {A_q}{(1-2q)} \approx \frac {(1-2r)R_q}{(1-2q)R_r}=c_{q,r}\pi _{q,r}~,
\end{equation}
 where $c_{q,r}=(1-2r)/(1-2q)<1.$ Hence  $\pi _{q,r}\approx \hat \pi _{q,r}\equiv  A_q/(1-2r)$ is approximately monotone increasing in the peakedness measure $\eta  _q,$ justifying the name \lq measure of peakedness\rq.\ Substitution of $c_{q,r}\hat \pi _{q,r}$ for $A_q/(1-2q)$ in (\ref{hornext}) yields an approximation
for $\eta  _{q,r}$ that is hereafter denoted $\hat \eta  _{q,r}$.  Examples are shown as thick dashed lines in Figure~\ref{fig2}.

\subsection{Tail-weight Measures}\label{tailweight}
To justify calling $\tau_{p,q}$ a tail-weight measure, recall that $F$ has a right tail with (asymptotic)
index $\alpha _R>0$ if $1-F(x)\sim u(x)x^{-\alpha _R}$ as $x\to \infty,$ where $u(x)$ is a slowly varying function.   Noting that density estimation for wide-tailed distributions
is difficult,  \cite{M-T-2000} introduce what
they call a \lq practical tail index\rq ,\ which, in our notation, for  $0<p<q<0.5$
is the ratio $\alpha _R(p,q)=\ln (q/p)/\ln(x_{1-p}/x_{1-q}).$  They explain why this gives a good indication of the size of $\alpha _R$, especially when computed for a range of pairs $(p,q).$ Similarly, if the left
tail index is denoted $\alpha _L$, one can derive $\alpha _L(p,q)=\ln (q/p)/\ln(x_{p}/x_{q}).$  It follows that $x_p=x_q\,(q/p)^{1/\alpha_L(p,q)}$ and $x_{1-p}=x_{1-q}\,(q/p)^{1/\alpha_R(p,q)}$, so that
\begin{equation}\label{tauMT}
    \tau _{p,q}=\frac {R_p}{R_q}=\frac {x_{1-q}\,\left (\frac {q}{p}
\right )^{1/\alpha _R(p,q)}-x_q\,\left (\frac {q}{p}\right )^{1/\alpha _L(p,q)}}{x_{1-q}-x_q} ~.
\end{equation}
This expression shows how the left and right hand practical tail indices  affect $\tau _{p,q}.$ For example as $\alpha _R(p,q)$ grows large, indicating a short right tail, the first term in the  numerator           of (\ref{tauMT}) approaches $x_{1-q}$\,, the first term in the denominator. But if $\alpha _R(p,q)$ decreases, the same first term of the numerator becomes larger than the first term below it. Similar remarks can be made for the left tail, but the main point is the $\tau _{p,q}$ increases as either of the practical tail indices decrease, as one would expect of a measure of tail-weight.

 For symmetric distributions, $\alpha _L(p,q)=\alpha _R(p,q)\equiv \alpha (p,q)$, so $\tau _{p,q}=(q/p)^{1/\alpha (p,q)}.$
\cite{M-T-2000} give
examples, including the Student-$t_\nu $ distribution which has tail index $\alpha =\nu ,$ and their
$H_h$ distributions for which $\alpha =1/h.$
For such distributions moments of larger order than $\alpha $ do not exist.

\subsection{Examples of Peakedness and Tail-weight}\label{quartiles}

A distribution-free choice for partitioning the distribution
is $x_{0.125}$, $x_{0.25}$ and $x_{0.375}$, so that comparisons are made between the ranges of the central quarter, half and three-quarters of the population.
 Table~\ref{table2} tabulates values of $\pi_{q,r},$  $\tau_{p,q}$ and $\kappa _{p,r}$ for this partition;
that is, for $p=1/8$, $q=1/4$ and $r=3/8.$  The symmetric models are listed in terms of increasing values
of $\kappa _{p,r}$ and similarly for the asymmetric models. Note that peakedness $\pi _{q,r}$ contributes more than tail-weight $\tau _{p,q}$ for all models except Models~10, 20 and 21. Models 10, the Cauchy, and the
skewed Cauchy $t_{2,1}$ have identical values and peakedness and tail-weight contribute equally to kurtosis
for each of them, as guaranteed by the results in \cite{JRP-2011}. Only Model 21, the very skewed $t_{2,1/2}$ family, has a larger tail-weight than peakedness.

Perhaps it is worth noting that the peakedness of $t_\nu $ and $\chi ^2_\nu $ models increases
with decreasing $\nu $, as one would expect from comparison of graphs of their densities.
 The only drawback of these definitions in terms of the ranges of the middle quarter, half and three-quarters
of the population is that the kurtosis for the normal model does not agree with the classical measure; here
the uniform model has kurtosis equal to 3.  Also, for the normal model there is not much difference between
the peakedness and tail-weight, and traditionalists might expect that tail-weight should contribute much less
than peakedness, because the normal model has relatively short tails.

\begin{table}[t!]
\caption{\label{table2}\em \footnotesize Columns 2--4 give the quantile peakedness $\pi_{q,r}$, the quantile tail-weight  $\tau_{p,q}$ and
their product, the kurtosis $\kappa _{p,r}$, for various models $F$ when $p=1/8$, $q=1/4$ and $r=3/8.$
Columns 5--7 contain the corresponding values when $r=1/3,$ $p=\Phi(3\Phi^{-1}(r))\approx
 0.1$,  $q=\Phi ^{-1}(-1)\approx 0.158$. }
\begin{footnotesize}
\begin{center}
\begin{tabular}{lrrrrrrr}
\hline
$\qquad F$   &  $\qquad \pi_{q,r}\ $ &  $\ \tau_{p,q}\ $ & $\quad  \kappa _{p,r}\ $
& $\quad $ & $\ \pi_{q,r}\ $ &  $\ \tau_{p,q}\ $ &  $\quad  \kappa _{p,r}\ $  \\
  \hline
1. Beta$(1/2,1/2)$    &  1.848 & 1.307 &  2.414 & & 1.757 & 1.085 & 1.906  \\
2. Uniform            &  2.000 & 1.500 &  3.000 & & 2.048 & 1.177 & 2.411  \\
3. Beta$(2,2)$        &  2.064 & 1.606 &  3.316 & & 2.193 & 1.235 & 2.709  \\
4. Normal             &  2.117 & 1.706 &  3.610 & & 2.322 & 1.292 & 3.000  \\
5. Logistic           &  2.151 & 1.771 &  3.809 & & 2.407 & 1.330 & 3.200  \\
6. Student-$t_5$      &  2.158 & 1.790 &  3.864 & & 2.429 & 1.342 & 3.260  \\
7. Student-$t_4$      &  2.170 & 1.815 &  3.938 & & 2.460 & 1.357 & 3.337  \\
8. Student-$t_2$      &  2.236 & 1.964 &  4.392 & & 2.643 & 1.446 & 3.820  \\
9. Laplace            &  2.409 & 2.000 &  4.819 & & 2.831 & 1.418 & 4.015  \\
10. Cauchy            &  2.414 & 2.414 &  5.828 & & 3.182 & 1.709 & 5.438  \\
 \noalign{\smallskip}
11. Beta$(2,1)$       &  2.054 & 1.590 &  3.265 & & 2.170 & 1.226 &  2.661\\
12. $\chi ^2_5$       &  2.127 & 1.725 &  3.669 & & 2.347 & 1.304 &  3.060\\
13. $\chi ^2_3$       &  2.136 & 1.743 &  3.722 & & 2.370 & 1.314 &  3.113\\
14. $\chi ^2_2$       &  2.151 & 1.771 &  3.809 & & 2.407 & 1.330 &  3.200\\
15. $\chi ^2_1$       &  2.229 & 1.906 &  4.249 & & 2.595 & 1.397 &  3.625\\
16. Log-normal        &  2.243 & 1.956 &  4.386 & & 2.646 & 1.432 &  3.789\\
17. Skew-$t_{2,2}$    &  2.236 & 1.964 &  4.392 & & 2.643 & 1.446 &  3.820\\
18. Pareto(2)         &  2.296 & 2.081 &  4.780 & & 2.800 & 1.506 &  4.216\\
19. Skew-$t_{2,1}$    &  2.414 & 2.414 &  5.828 & & 3.182 & 1.709 &  5.438\\
20. Skew-$t_{2,1/2}$  &  2.996 & 4.222 & 12.649 & & 5.329 & 2.633 & 14.033\\
   \end{tabular}
  \end{center}
  \end{footnotesize}
     \end{table}

A Gaussian-centric choice could define the central portion of the distribution as that lying within one
standard deviation of the mean; that is,  $q=\Phi ^{-1}(-1)= 0.1586553.$ Then, taking $r=1/3$ and
$p=\Phi(3\Phi ^{-1}(r))=0.098\approx 0.1$ for reasons given in Section~\ref{intro} gives somewhat different
results, also listed in Table~\ref{table2}.  Now the kurtosis for the normal is 3 by definition, and
the contribution of its peakedness factor is almost twice that of tail-weight.  In fact the contribution of
peakedness to tail-weight has increased for all distributions. Nevertheless, the orderings of kurtosis
within symmetric and asymmetric groups remains unchanged from the \lq model-free\rq\ choice of $p,$ $q$ and $r$.

\section{DISTRIBUTION-FREE INFERENCE}\label{Funknown}

The material in this section focusses on the kurtosis coefficient $\kappa _{p,r}$, but equally applies
to peakedness $\pi _{q,r}$ or tail-weight $\tau _{p,q}.$  Let $X_{([nr])}$ denote the $[nr]$th order statistic of a sample of size $n$ from $F$, and define the sample version of $R_r$ by $\hat R_r=R_r(F_n)=X_{(n-[nr]+1)}-X_{([nr])}.$
  We estimate $\kappa _{p,r}=R_p/R_r$ by $\hat\kappa _{p,r}=\hat R_p/\hat R_r.$

\subsection{Variance Stabilization of $\hat \kappa _{p,r} $}\label{VST}

The methodology for finding a variance stabilizing transformation (VST) of a ratio of statistics, each of which is a finite linear combination of order statistics, has already been established for other ratios of linear combinations of quantiles in \cite{S-2013a,S-2013b,S-2014}, so here we only restate the required results. One first shows
 that  $\var [\hat \kappa _{p,r}]=\var [\hat R_p/\hat R_r]$ satisfies
$ \var [\hat \kappa _{p,r}]   \doteq  \var [\hat R_p-\e [\hat \kappa _{p,r}] \,\hat R_r]/R_r^2.$
Therefore  $n\var [\hat \kappa _{p,r}]=q(\e [\hat \kappa _{p,r}] )$ where  $q(t) = a_0 +a_1t+a_2t^2$ is a  quadratic with constants:
\begin{eqnarray}\label{constants}
\nonumber
 a_0  &=& a_0(p,r)\;=\; n\var _F [\hat R_p]/R_r^2 \\
 a_1  &=&  a_1(p,r)\;=\; -2n\cov _F[\hat R_p,\hat R_r]/R_r^2 \\ \nonumber
 a_2  &=&\quad a_2(r)\;=\;  n \var _F[\hat R_r]/R_r^2 ~.
\end{eqnarray}
Note that $a_0,$ $a_1$ and $a_2$ are free of location, scale and sample size. The quadratic
$q(t)>0$ for all $t$  because $a_0>0$ and its discriminant $a_1^2-4a_0a_2<0$; the latter inequality follows from  $|a_1/\{2\sqrt{a_0a_2}\,\}|=|\corr [\hat R_p,\hat R_r]|< 1.$ Hereafter let $D^2=4a_0a_2-a_1^2$.
In the remainder of this subsection, drop the subscripts $p,r$ on $\kappa _{p,r}$.
A variance stabilizing transformation (VST) of $\hat \kappa $ is
\begin{equation}\label{kfun}
    h_n(x)= \sqrt {\frac{n}{a_2}}\;  \sinh ^{-1}\left \{\frac {q\,'(x)}{D}\right \} +c~,
\end{equation}
where $c$ is an arbitrary real number. In carrying out inference for $\kappa  $,
it is useful to center $h_n(\hat \kappa )$ at an arbitrary null hypothesis value $\kappa _0\geq 1$
by introducing  $T_{n,\kappa _0}(\hat \kappa )=h_n(\hat \kappa )-h_n(\kappa _0)$, so $\e _{\kappa _0}[T_{n,\kappa _0}]$
is approximately 0 under the null.  By the Delta Theorem \cite[p.40]{Das-2008}, as $n$ grows without bound,
$T_{n,\kappa _0}(\hat \kappa )\sim N(\sqrt n\;K _{\kappa _0}(\kappa ),1),$ where
\begin{equation}\label{key}
    K _{\kappa _0}(\kappa )= \frac {1}{\sqrt {a_2}\,}\left [\sinh ^{-1}\left \{\frac {q\,'(\kappa )}{D}\right \}- \sinh ^{-1}\left \{\frac {q\,'(\kappa _0)}{D}\right \}\right ]~.
\end{equation}
We can write $T_{n,\kappa _0}(\hat \kappa )= \sqrt n\;K _{\kappa _0}(\hat \kappa ).$  A level-$\alpha $ test rejects the null $\kappa =\kappa _0$ in favor of $\kappa >\kappa _0$ for $T_{n,\kappa _0}(\hat \kappa )\geq z_{1-\alpha }=\Phi ^{-1}(1-\alpha ).$

To make this statistic distribution-free, the nuisance parameters $a_0$, $a_1$ and $a_2$ must be estimated. They depend on the unknown $F$ through the sparsity index $g_p=g(p)
=1/f(x_p) $ of \cite{tukey-1965}, at each of the quantiles $x_p<x_r<x_{1-r}<x_{1-p}$\,. This requires density estimates at the selected
quantiles, and the resulting constants are denoted  $\hat a_0$, $\hat a_1$ and $\hat a_2$. When these estimated constants are substituted into $q$, $D$, and $T_{n,\kappa _0}(\hat \kappa )= \sqrt n\;K _{\kappa _0}(\hat \kappa )$, the results are denoted $\hat q$, $\hat D$ and $T_{\text{DF},n,\kappa _0}(\hat \kappa )$.
The method of sparsity density estimation described in \cite[Sec. 4.1]{S-2014} is also utilized here;
but the constants (\ref{kurtconst}) are different in this kurtosis setting.

\subsection{Constants required by the VST}\label{vstconstants}

For fixed $0< r\leq s <1$ and sample size $n$ increasing without bound,  $\e [X_{([nr]}]  \doteq x_r $ and $
  n\cov  [X_{([nr]},X_{([ns]}] \doteq r(1-s)g_rg_s~,$ where \lq$ \doteq $\rq \  means that lower order terms
are ignored; see, eg. \cite[p.80]{david-1981} or \cite[p.93]{Das-2008}.

It follows that for $0<p<r<1/2$ the constants (\ref{constants}) required by the VST are:
\begin{eqnarray}\label{kurtconst}
\nonumber
  R _r^2\,a_0(p,r) &=&  p(g_p^2+g_{1-p}^2)-p^2(g_p+g_{1-p})^2\\ \nonumber
  R _r^2\,a_1(p,r) &=&  2\{pr(g_r\,g_{1-p}+g_p\,g_{1-r})-p(1-r)(g_p\,g_r+g_{1-p}\,g_{1-r}) \}\\
  R _r^2\,a_2(r)   &=&  r(g_r^2+g_{1-r}^2)-r^2(g_r+g_{1-r})^2~.
\end{eqnarray}
When $F$ is symmetric, $R_r=2x_{1-r}$ and $g_r=g_{1-r}$, so these formulae reduce to $a_0(p,r)=2p\,g_p^2/R _r^2,$
$a_1(p,r)= 4p\,g_p\,g_r\,(2r-1)/  R _r^2$ and $a_2(r)=2r\,g_r^2/  R _r^2.$
Table~\ref{table3} lists values of  $\kappa_{1/3}=\kappa _{p(1/3),1/3}$, where $p(r)=\Phi(3\Phi^{-1}(r)),$ and the VST constants  $a_0$, $a_1$ and $a_2$.

\subsection{Two-sided Confidence Intervals for $\kappa $}\label{sec:DFcis}

A nominal 100$(1-\alpha)$\% distribution-free confidence interval for   $\kappa $ is derived exactly as  for the skewness coefficient in \cite[Sec. 3.3]{S-2014} and displayed in Equation~9 of that paper; its
analogue here is, for $c_\alpha =z_{1-\alpha /2}$\, :
\begin{equation}\label{DFcis}
    [L,U]_\text{DF}=\frac {1}{\hat a_2}\left [ \hat D\sinh\left\{\sinh^{-1}\left(\frac {\hat q\,' (\hat \kappa )}{\hat D}\right )
    \pm c_\alpha \sqrt {\frac {\hat a_2}{n}}\, \right \}-\hat a_1\right ]~.
\end{equation}
In this expression $\hat a_0$, $\hat a_1$ and $\hat a_2$ as well as $\hat q^{\,\prime}$ and $\hat D$ are all estimated
using distribution-free methods.
The empirical coverage of nominal 90\% and 95\% distribution-free
confidence intervals for $\kappa $ based on (\ref{DFcis})
are found for various $n$ in Section~\ref{simulations}.
Also of interest are the widths of these intervals, defined by $W=U_{DF}-L_{DF} $.

\begin{table}[t!]
\caption{\label{table3}\em \footnotesize For $r=1/3$ and $p= \Phi (3\Phi ^{-1}(r))=0.1$ are listed the kurtosis coefficient $\kappa _{1/3}=R_p/R_{r}$, the VST constants (\ref{constants}), the asymptotic width  $w_{asym}=2\sqrt {q (\kappa _{1/3})}\, $ appearing in (\ref{What}) and the asymptotic relative widths $ rw_{asym}
= w_{asym}/\kappa _{1/3}.$}
\begin{small}
\begin{center}
\begin{tabular}{lrrrrrr}
\hline
  $\qquad F$ &  $\kappa_{1/3}$  & $a_0\ $ & $a_1\quad $ & $\quad a_2\quad $ & $w_{asym}$ &  $rw_{asym}$\\
\hline
1. Beta$(1/2,1/2)$ &   1.906 &   0.143 & $ -0.339$ & 1.645 &  4.678  & 2.455 \\
2. Uniform         &   2.411 &   1.420 & $ -1.178$ & 2.000 &  6.390  & 2.650 \\
3. Beta$(2,2)$     &   2.709 &   3.512 & $ -1.919$ & 2.146 &  7.499  & 2.770 \\
4. Normal          &   3.000 &   7.094 & $ -2.802$ & 2.265 &  8.735  & 2.912 \\
5. Logistic        &   3.200 &  10.478 & $ -3.462$ & 2.342 &  9.670  & 3.022 \\
6. Student-$t_5$   &   3.260 &  11.882 & $ -3.699$ & 2.358 &  9.975  & 3.060 \\
7. Student-$t_4$   &   3.337 &  13.646 & $ -3.986$ & 2.384 & 10.371  & 3.108 \\
8. Student-$t_2$   &   3.820 &  28.436 & $ -5.930$ & 2.531 & 13.073  & 3.422 \\
9. Laplace         &   3.200 &  20.049 & $ -5.772$ & 3.752 & 12.648  & 3.953 \\
10. Cauchy         &   5.438 & 137.680 & $-14.024$ & 2.924 & 24.323  & 4.473 \\
\noalign{\smallskip}
11. Beta$(2,1)$     &   2.661 &    4.088 & $  -2.261$ &  2.311 &   7.599 &  2.856 \\
12. $\chi ^2_5$     &   3.060 &   10.939 & $  -3.931$ &  2.543 &   9.532 &  3.115 \\
13. $\chi ^2_3$     &   3.113 &   14.104 & $  -4.857$ &  2.773 &  10.171 &  3.267 \\
14. $\chi ^2_2$     &   3.200 &   18.899 & $  -6.244$ &  3.122 &  11.115 &  3.474 \\
15. $\chi ^2_1$     &   3.625 &   41.492 & $ -12.353$ &  4.660 &  15.226 &  4.200 \\
16. Log-normal      &   3.789 &   49.077 & $ -11.552$ &  3.811 &  15.495 &  4.089 \\
17. Skew-$t_{2,2}$  &   3.820 &   54.245 & $ -12.480$ &  3.943 &  16.014 &  4.192 \\
18. Pareto(2)       &   4.216 &   90.352 & $ -17.572$ &  4.496 &  19.616 &  4.652 \\
19. Skew-$t_{2,1}$  &   5.438 &  282.221 & $ -32.651$ &  4.963 &  31.712 &  5.831 \\
20. Skew-$t_{2,1/2}$&  14.033 & 6958.645 & $-218.133$ &  8.456 & 149.167 & 10.630 \\
\end{tabular}
\end{center}
\end{small}
\end{table}

\clearpage
\newpage

They can be expressed, see \cite[App.2]{S-2014},
\begin{equation}\label{What}
 W= \frac {w_{asym}(\hat \kappa )\;z_{1-\alpha /2}}{\sqrt n\,}+o_p(n^{-1/2})~,
\end{equation}
where $w_{asym}(\kappa )=2\sqrt {q (\kappa )}$\; and $q(t) = a_0 +a_1t +a_2t^2$. Thus for large $n$ the  half-width of the confidence intervals (\ref{DFcis}) is approximately $z_{1-\alpha /2}$ times
the standard error  of $\hat \kappa $, which is $\sqrt {\var [\hat \kappa]}\approx \sqrt {q(\kappa )/n}\ $.

Since $\hat \kappa $ is consistent for $\kappa $ it is of interest to evaluate
 $2\sqrt {q (\kappa )}\, $ for various $F$.  It turns out that the interval widths are almost linearly increasing  with
$\kappa ,$ so we also introduce the relative width $rW=W/\kappa $.
It follows from (\ref{What}) that to obtain a large sample
100($1-\alpha )$\% confidence interval for $\kappa $ of desired relative width $rW_0=W_0/\kappa $ one requires
\begin{equation}\label{nmin}
n\geq n_0=n_0(\alpha ,rW_0)= \biggl \{\frac {\max _F \{rw_{asym}(F)\}\;z_{1-\alpha /2}}{rW_0}\biggr \}^2~.
\end{equation}
where $rw_{asym}(F)=2\sqrt {q (\kappa (F))}\,/\kappa (F).$ By referring to Table~\ref{table3}, one sees
that for the choices $r=1/3$, $q=1/10$, excluding the skew-$t$ distributions, $rw_{asym}(F)\leq 4.652$. To ensure $rW_0=0.2$ with 95\% confidence, one requires $n\geq n_0=(4.652\times 1.96\times 5)^2=2079$.

\section{SIMULATION RESULTS}\label{simulations}

\subsection{Empirical Coverage and Widths}

We find distribution-free confidence intervals for $\kappa =\kappa _{p(r),r}$, where $r=1/3$ and
$p(r)=\Phi(3\Phi ^{-1}(r))\approx 0.1,$ for reasons given in Section~\ref{intro}.
In our simulation studies we used the software package R \cite{R}, and estimated the sparsity index using
the method described in \cite[Sec.4]{S-2014}.
 \begin{table}[t!]
\caption{\label{table4}\em \footnotesize  Estimates of coverage probabilities and widths of nominal 90\% and 95\% confidence intervals for $\kappa _r=\kappa _{p(r),r}$ when $r=1/3$ and $p(r)=\Phi(3\Phi ^{-1}(r))= 0.09815$, all based on 40,000 replications of samples from selected symmetric models. The average interval relative widths are not shown but can be recovered from
$\overline {rW}=\widehat {rw}\;z_{1-\alpha /2}/\sqrt{n}\,,$ see (\ref{What}).}
\begin{footnotesize}
\begin{center}
\begin{tabular}{lrrrrcrrrr}
  &\multicolumn{4}{c}{\rule[-1mm]{0mm}{6mm}90\%} &  &\multicolumn{4}{c}{95\%}\\
\hline
$F$ &  $n$ &  $\hat {\kappa}_{1/3}$&  $cp\quad $  &$\widehat {rw}$ & & $n$ & $\hat {\kappa}_{1/3}$&  $cp\quad $ & $\widehat {rw}$ \\
\noalign{\smallskip}\hline\noalign{\smallskip}
               &  100  & 2.454  & 0.918 & 2.874  &   &   100  &  2.454 & 0.961 &  2.887 \\
2. Uniform     &  400  & 2.421  & 0.906 & 2.715  &   &   400  &  2.421 & 0.954 &  2.719 \\
               & 1000  & 2.415  & 0.904 & 2.679  &   &  1000  &  2.415 & 0.953 &  2.680 \\
               & 4000  & 2.412  & 0.901 & 2.658  &   &  4000  &  2.412 & 0.950 &  2.659 \\
             &$+\infty $& 2.411 & 0.900 & 2.650  & & $+\infty $& 2.411 & 0.950  &  2.650  \\
\noalign{\medskip}
               & 100  & 3.060  & 0.928 & 3.227  &   &    100  &  3.057 & 0.966 &  3.239 \\
4. Normal      & 400  & 3.015  & 0.913 & 3.020  &   &    400  &  3.014 & 0.957 &  3.022 \\
               &1000  & 3.005  & 0.906 & 2.965  &   &   1000  &  3.007 & 0.953 &  2.966 \\
               &4000  & 3.001  & 0.903 & 2.933  &   &   4000  &  3.001 & 0.954 &  2.933 \\
               &$+\infty $&  3.000 & 0.900 &  2.912 &   & $+\infty $&  3.000 & 0.950 &  2.912  \\
\noalign{\medskip}
               &  100   & 3.334  & 0.930  & 3.382  &   &   100  &  3.337 & 0.967 &  3.400 \\
6. Student-$t_5$&  400  & 3.276  & 0.909  & 3.160  &   &   400  &  3.275 & 0.957 &  3.164 \\
               & 1000   & 3.267  & 0.904  & 3.106  &   &  1000  &  3.266 & 0.953 &  3.107 \\
               & 4000   & 3.261  & 0.901  & 3.077  &   &  4000  &  3.261 & 0.951 &  3.077 \\
               & $+\infty $&  3.260 & 0.900 & 3.060 &   & $+\infty $& 3.260 & 0.950 &  3.060 \\
\noalign{\medskip}
               &  100   & 3.949  & 0.927  & 3.772  &   &   100  & 3.945  & 0.965 & 3.791   \\
8. Student-$t_2$  &  400   & 3.850  & 0.906  & 3.503  &   &   400  & 3.853  & 0.954 & 3.507   \\
               & 1000   & 3.831  & 0.903  & 3.446  &   &  1000  & 3.834  & 0.953 & 3.449   \\
               & 4000   & 3.823  & 0.901  & 3.423  &   &  4000  & 3.822  & 0.951 & 3.424   \\
               & $+\infty $&  3.820 & 0.900 & 3.422 &   & $+\infty $&  3.820 & 0.950 &  3.422  \\
\noalign{\medskip}
                 &    100   & 5.806  & 0.905  & 4.996   &   &   100  &  5.797  & 0.948 & 5.031  \\
10. Cauchy         &    400   & 5.523  & 0.900  & 4.585   &   &   400  &  5.525  & 0.948 & 4.595  \\
               &   1000   & 5.473  & 0.899  & 4.488   &   &  1000  &  5.471  & 0.948 & 4.492  \\
               &   4000   & 5.446  & 0.900  & 4.459   &   &  4000  &  5.448  & 0.948 & 4.459  \\
          &   $+\infty $  &  5.438 & 0.900  & 4.473   &   & $+\infty $& 5.438 & 0.950 & 4.473 \\
\end{tabular}
\end{center}
\end{footnotesize}
\end{table}
\clearpage
\newpage

 \begin{table}[t!]
\caption{\label{table5}\em \footnotesize  Estimates of coverage probabilities and widths of nominal 90\% and 95\% distribution-free confidence intervals for $\kappa _r=\kappa _{p(r),r}$ when $r=1/3$ and $p(r)=\Phi(3\Phi ^{-1}(r))\approx 0.1$, for selected asymmetric models. Notation as in Table~\ref{table4}.}
\begin{footnotesize}
\begin{center}
\begin{tabular}{lrrrrcrrrr}
  &\multicolumn{4}{c}{\rule[-1mm]{0mm}{6mm}90\%} &  &\multicolumn{4}{c}{95\%}\\
\hline
$F$ &  $n$ &  $\hat {\kappa}_{1/3}$&  $cp\quad $ & $\widehat {rw}$ & & $n$ & $\hat {\kappa}_{1/3}$&  $cp\quad $ & $\widehat {rw}$ \\
\noalign{\smallskip}\hline\noalign{\smallskip}
               &  100  & 3.124  & 0.924 & 3.380  &   &   100   &  3.127 & 0.963 & 3.398 \\
12. $\chi ^2_5$    &  400  & 3.076  & 0.908 & 3.189  &   &   400   &  3.076 & 0.955 & 3.194 \\
               & 1000  & 3.066  & 0.902 & 3.147  &   &  1000   &  3.066 & 0.954 & 3.148 \\
               & 4000  & 3.060  & 0.902 & 3.125  &   &  4000   &  3.061 & 0.951 & 3.125 \\
           & $+\infty $& 3.060  & 0.900 & 3.115  &   &$+\infty $& 3.060 & 0.950 & 3.115 \\
\noalign{\medskip}
               & 100  & 3.279  & 0.911 & 3.654   &   &    100  & 3.283 & 0.966 & 3.676  \\
14. $\chi ^2_2$    & 400  & 3.217  & 0.901 & 3.494   &   &    400  & 3.220 & 0.957 & 3.499  \\
               &1000  & 3.206  & 0.897 & 3.461   &   &   1000  & 3.209 & 0.953 & 3.468  \\
               &4000  & 3.202  & 0.899 & 3.457   &   &   4000  & 3.202 & 0.954 & 3.458  \\
           &$+\infty $&  3.200 & 0.900 & 3.474   &   & $+\infty $&  3.200 & 0.950 & 3.474  \\
\noalign{\medskip}
               &  100   & 3.912  & 0.900 & 4.230  &   &   100  &  3.911 & 0.945 &  4.249 \\
16. Lognormal     &  400   & 3.821  & 0.893 & 4.083  &   &   400  &  3.820 & 0.945 &  4.084 \\
               & 1000   & 3.802  & 0.895 & 4.052  &   &  1000  &  3.801 & 0.948 &  4.059 \\
               & 4000   & 3.792  & 0.897 & 4.059  &   &  4000  &  3.792 & 0.950 &  4.059 \\
               & $+\infty $&  3.789  & 0.900 & 4.089 &   & $+\infty $& 3.789 & 0.950 &  4.089 \\
\noalign{\medskip}
               &  100   & 4.394  & 0.927  & 4.729 &   &   100  & 4.397 & 0.932 &  4.772  \\
18. Pareto(2)      &  400   & 4.264  & 0.906  & 4.634 &   &   400  & 4.261 & 0.939 &  4.639  \\
               & 1000   & 4.232  & 0.903  & 4.599 &   &  1000  & 4.235 & 0.943 &  4.607  \\
               & 4000   & 4.220  & 0.901  & 4.611 &   &  4000  & 4.221 & 0.946 &  4.611  \\
               & $+\infty $&  4.216 & 0.900 & 4.652 &   & $+\infty $&  4.216 & 0.950 &  4.652  \\
\noalign{\medskip}
               &    100   & 5.831  & 0.847  &  $-$  &   &   100  &  5.817  &0.900 &  $-$  \\
19. Skew-$t_{2,1}$ &    400   & 5.527  & 0.868  & 5.852  &   &   400  &  5.525  &0.922 & 5.874  \\
               &   1000   & 5.471  & 0.883  & 5.775  &   &  1000  &  5.473  &0.932 & 5.770  \\
               &   4000   & 5.447  & 0.893  & 5.794  &   &  4000  &  5.448  &0.943& 5.794  \\
               &$+\infty $& 5.438  & 0.900  & 5.831  &   & $+\infty $&  5.438 & 0.950 &  5.831 \\
\end{tabular}
\end{center}
\end{footnotesize}
\end{table}

\clearpage
\newpage

In Table~\ref{table4} are shown the results of 40,000 simulations from 5 symmetric models with sample sizes
ranging from 100 to 4000. For each replicate $\hat \kappa _{1/3}$ and the VST constants $a_0, a_1$ and $a_2$ of (\ref{constants}) were also estimated, and a confidence interval found using (\ref{DFcis})
with these estimated constants. The average value of these estimates $\hat \kappa _{1/3}$  is shown in Column~3 of Table~\ref{table4}. Note the positive bias for smaller $n$. Despite this bias, the
empirical coverage probabilities of $\kappa _{1/3}=2.411$ in Column~4 are only slightly conservative for $n\geq 100.$ To obtain the estimates of $rw_{asym}$ shown in Column~5, we found the average of
    $\widehat {rw}_{asym}=\sqrt{n}\,(rW)/z_{0.95}$, where $rW=(U-L)/\hat\kappa_{1/3},$
    see formulae (\ref{What}). Note
that these estimates are also converging  to their limiting value, shown in the last row for each model, and obtained from Table~\ref{table3}. Similar results are found for 95\% confidence intervals,  listed in the right hand columns of Table~\ref{table4}.

For simulated data generated from the normal
and Student-$t_5$ models with  sample sizes 100 the coverage probabilities are conservative,
but for 400 or more the coverages and widths are reflecting our expectations.
For heavier tailed distributions such as the Student-$t_2$ and Cauchy distributions,
 the methods can fail to work at all for the smaller sample sizes. This is because
outliers  in the samples can undermine the estimates of the sparsity index, occasionally leading to negative
values of $\widehat {D^2}=4\hat a_0\hat a_2-\hat a_1^2.$

 R software functions for
finding the VST constants and the resulting distribution-free confidence intervals are available online, see Section~\ref{sec:supp}.

In Table~\ref{table5} the results of similar studies for five asymmetric models are presented.
Again, $\hat \kappa _{1/3}$ is biased upwards, but converges to its target $\kappa _{1/3}$. Now the
 sample sizes $n\geq 400$ appear necessary to obtain 90 or 95\% confidence intervals, as the case may be.
For $n=100$ and Model~19, no empirical average widths are tabled for the reasons just given in the last paragraph.  For Model~19, which has the same kurtosis as Model 10, at least 40 times as many observations are required to obtain accurate coverage. Thus estimating $\kappa _{1/3}$ can be costly for very skewed distributions.

\begin{figure}[t!]
\centering
\vspace{-1cm}
\includegraphics[scale=.6]{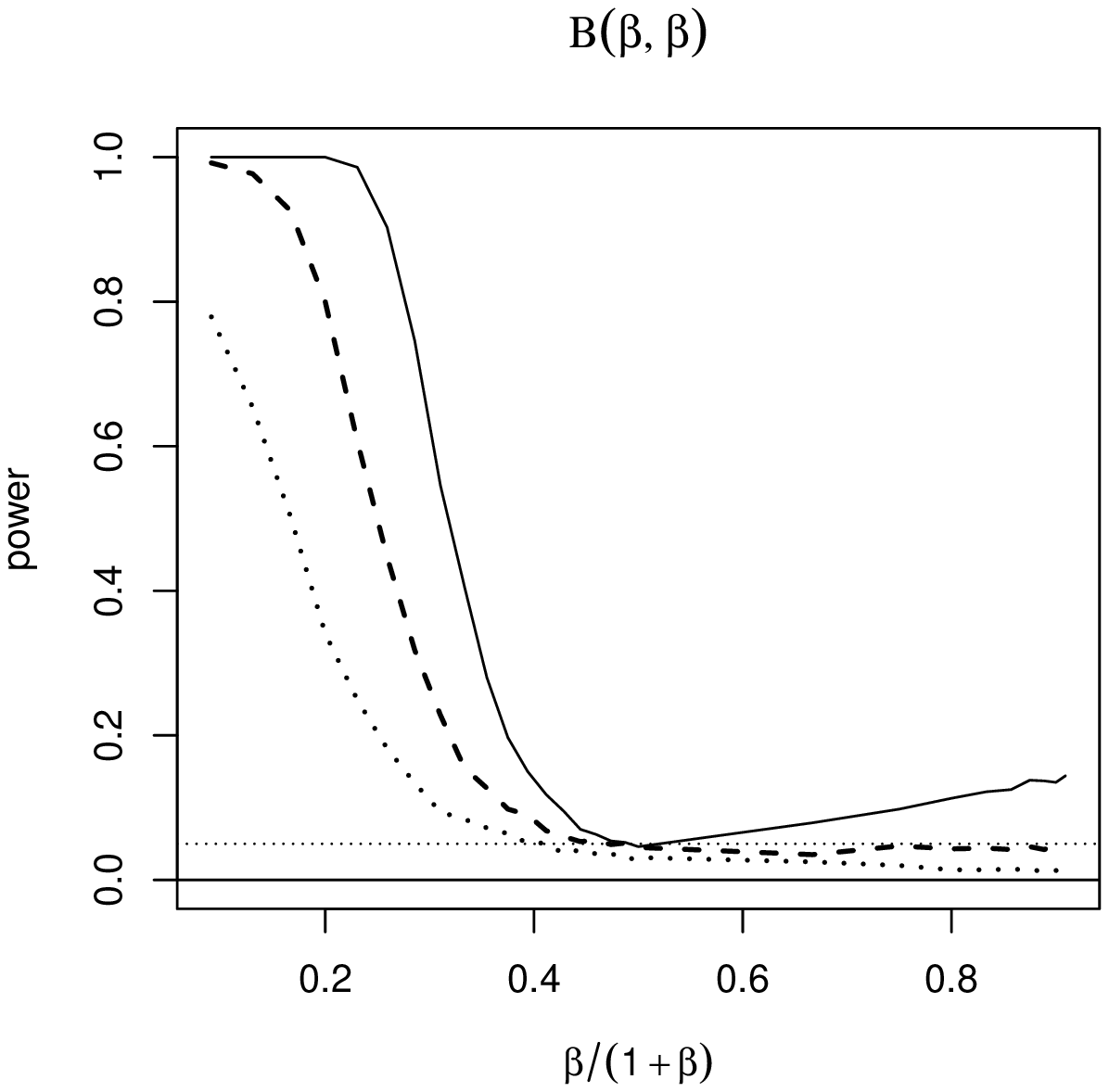}
\vspace{-0.5cm}
\caption{\footnotesize  Graphs of empirical power of two-sided level-0.05 tests for the $\Beta (\beta ,\beta )$ model
plotted as a function of $\beta/(\beta +1).$  The power for $n=50$
is shown as a dotted line, for $n=200$ as a dashed line and for $n=800$ as a solid line. The dotted horizontal line gives
the level of the test. \label{fig3}}
\includegraphics[scale=.6]{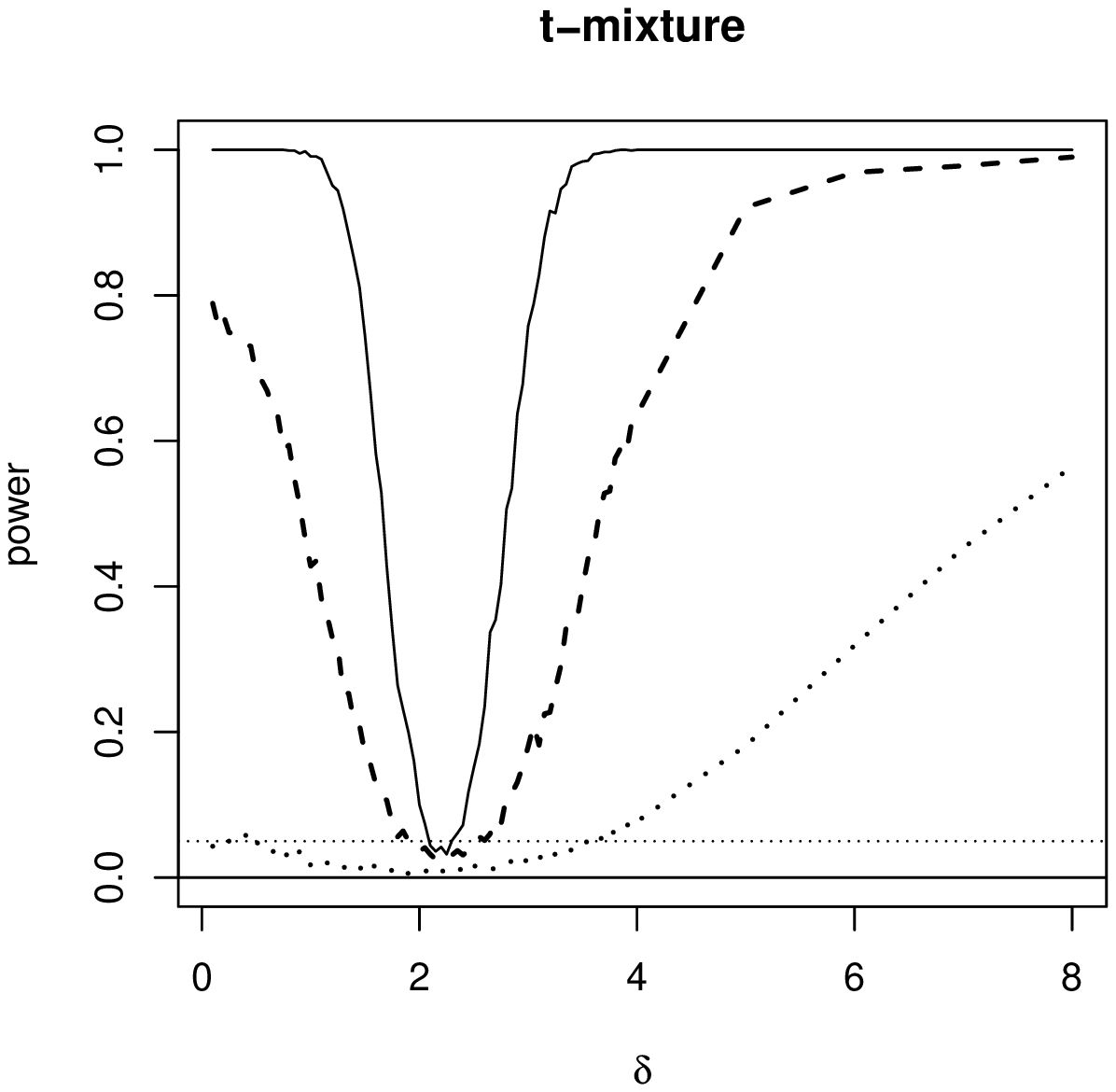}
\vspace{-0.5cm} \caption{\footnotesize  As in Figure~\ref{fig3}, but now for a 50:50 mixture of two $t_{1/2}$ distributions that are distance $\delta $ apart.
\label{fig4}}
\end{figure}

\subsection{Power of $\hat \pi _{q,r}$ for  Detecting Bimodality}

In this section we illustrate the power of $\hat \pi _{q,r}$ to detect bimodality, as well as peakedness.
  In Section~\ref{peakedness} it was shown that the extended Horn's peakedness
measure  defined in (\ref{hornext}) is monotone in $\pi _{q,r}$ via (\ref{ratioapprox}).
Recall that $\eta _q$ lies in  $[-1,1]$ with negative values indicating \lq bimodality', positive values
\lq peakedness' and 0 \lq uniformity' near the median. Therefore a two-sided test of $\eta _q=0$  is approximately a two-sided test of $\pi _{q,r}=(1-2q)/(1-2r).$

Fix $q=1/4$, $r=3/8$.    In Figure~\ref{fig3} is shown the empirical power of the level-0.05 test of $\pi _{1/4,3/8}=2,$ when the data
are generated from the symmetric $\Beta (\beta ,\beta )$ model, for selected values of $\beta $. The power for $n=50$
is shown as a dotted line, for $n=200$ as a dashed line and for $n=800$ as a solid line. These curves are near 0.05
when $\beta =1$, (the uniform model). For example, the power of detecting the bimodal model $\Beta (1/3,1/3)$ is approximately 0.4 for $n=200$ and 0.8 for $n=800 $ observations. There is not much power for detecting peakedness for large $\beta $
because the $\Beta (\beta ,\beta )$ model approaches the Normal as $\beta \to \infty .$

\clearpage
\newpage

Figure~\ref{fig4} shows  the empirical power of the same distribution-free test for detecting peakedness
and bimodality of a 50:50 mixture of two central $t_{1/2}$ distributions, as a function of the distance
$\delta $ between them.  For sample size $n=200$ the test has the right level and the power to detect either
peakedness or bimodality, as the case may be, depending on the value of $\delta .$

\section{FURTHER RESEARCH}\label{summary}

We extended the peakedness measure of \cite{horn-1983} to arbitrary densities and showed
that the ratio of interquantile ranges of  \cite{rupp-1987} is approximately monotone in it when
applied to the central portion of the distribution; that is, for $\pi _{q,r}$, where $q=\Phi^{-1}(-1)\approx 0.16$
or $q=0.25$, say, and $q<r<0.5$. When applied to  the non-central portion,  Ruppert's ratio $\tau _{p,q}$ for $0<p<q$
is also monotone in the practical tail index of \cite{M-T-2000}.
We endorse the idea that peakedness and tail-weight are best estimated separately, because as
the simple factorization $\kappa _{p,r}=\pi _{q,r}\tau _{p,q}$ shows, kurtosis is fundamentally
a product of peakedness and tail-weight.

Distribution-free confidence intervals are derived for $\kappa _{p,r}$, and hence
available for $\pi _{q,r}$ and $\tau _{p,q}$ separately. In our simulation studies
we concentrated on estimation of  $\kappa _{1/10,1/3}$, and shown that it is possible to obtain accurate 90\% and 95\% distribution-free confidence intervals for data simulated from a large variety of distributions, provided that the sample sizes were at least 400. This procedure is resistant to almost 10\% outliers on either side of the sample. A formula for choosing the sample size required to obtain intervals of a given desired relative width over a large class of models is included.

\cite{schmid-2003} found finite-sample and asymptotic tests for normality based on $\hat \pi _{1/8,1/4}$, $\hat \tau _{1/40,1/8}$ and $\hat \kappa _{1/40,1/4},$ using the asymptotic bivariate normality of the
sample interquantile ranges. The VST-transformed pair $(K_1,K_2)=(\key _{\pi _0}(\hat \pi _{q,r}),\, \key _{\tau _0}(\hat \tau _{p,q})),$ where $\key $ is of the form (\ref{key}), is also asymptotically
bivariate normal, with a covariance structure dependent on the sparsity indices at six quantiles.
Thus it should be possible to find 100$(1-\alpha )$\% distribution-free
confidence ellipses for the transformed
pair $(\key _{\pi _0}(\pi _{q,r}),\, \key _{\tau _0}(\tau _{p,q}))$, and, by back-transformation, non-elliptical
confidence regions for $(\pi _{q,r}, \tau _{p,q})$.

These methods can be used to find confidence intervals for the octile based kurtosis measure of \cite{moors-1988} and the quintile based kurtosis of \cite{JRP-2011}, and it would be of interest to
see whether or not they perform better than those presented here. A closely related problem is the
 estimation of tail-indices, and these methods can be easily adapted to find confidence intervals for the robust measures of
tail weights proposed by
\cite{brys-2006}. Extensions of these inferential methods to the multivariate setting are also of interest, see \cite{wang-2005}.

\section{SUPPLEMENTARY MATERIAL}\label{sec:supp}

Given a vector of data $x$, selected values $0<p <r <0.5$, and $\alpha $, this script will enable the user to find a 100$(1-\alpha )$\% distribution-free confidence interval for
Ruppert's measure of kurtosis $\kappa _{p,r}$, and hence also for the peakedness measure $\pi _{q,r}$ or the tail-weight measure $\tau _{p,q}.$

\begin{description}
  \item[findDFcikurt:\ ] R script
\end{description}

\clearpage
\newpage


\end{document}